\documentclass[mathpazo]{cicp}

\usepackage[utf8]{inputenc} 
\usepackage[T1]{fontenc}    
\usepackage{hyperref}       
\usepackage{url}            
\usepackage{booktabs}       
\usepackage{amsfonts}       
\usepackage{nicefrac}       
\usepackage{microtype}      
\usepackage{float}

\usepackage{amsmath,amssymb}
\usepackage{multirow}
\usepackage{diagbox}
\usepackage{graphicx}
\usepackage{subfigure}

\newcommand{\R}{{\mathbb{R}}}

\newcommand{\be}{\begin{eqnarray}}
\newcommand{\ee}{\end{eqnarray}}
\newcommand{\beq}{\begin{equation}\begin{aligned}}
\newcommand{\eeq}{\end{aligned}\end{equation}}
\newcommand{\beqn}{\begin{equation*}\begin{aligned}}
\newcommand{\eeqn}{\end{aligned}\end{equation*}}
\newcommand{\ben}{\begin{eqnarray*}}
\newcommand{\een}{\end{eqnarray*}}

\begin{document}
	
\title{A Neural Network with Plane Wave Activation for Helmholtz Equation}

%

\author[An Author et.~al]{Ziming Wang\affil{1,2}, 
	Tao Cui\affil{2, *}, and
	Xueshuang Xiang\affil{1,}\corrauth}

\address{\affilnum{1}\ Qian Xuesen Laboratory of Space Technology, 
	China Academy of Space Technology \\
	\affilnum{2}\ Academy of Mathematics and System Sciences, 
	Chinese Academy of Sciences}

\emails{{\tt tcui@lsec.cc.ac.cn} (T. Cui), {\tt xiangxueshuang@qxslab.cn} (X. Xiang)}

\begin{abstract}
This paper proposes a plane wave activation based neural network (PWNN) for solving Helmholtz equation, the basic partial differential equation to represent wave propagation, e.g. acoustic wave, electromagnetic wave, and seismic wave.  
Unlike using traditional activation based neural network (TANN) or $sin$ activation based neural network (SIREN) for solving general partial differential equations, 
we instead introduce a complex activation function $e^{\mathbf{i} {x}}$, the plane wave which is the basic component of the solution of Helmholtz equation. 
By a simple derivation, we further find that PWNN is actually a generalization of the plane wave partition of unity method (PWPUM) by additionally imposing a learned basis with both amplitude and direction to better characterize the potential solution. 
We firstly investigate our performance on a problem with the solution is an integral of the plane waves with all known directions. 
The experiments demonstrate that: PWNN works much better than TANN and SIREN on varying architectures or the number of training samples, 
that means the plane wave activation indeed helps to enhance the representation ability of neural network toward the solution of Helmholtz equation; 
PWNN has competitive performance than PWPUM, e.g. the same convergence order but less relative error. 
Furthermore, we focus a more practical problem, the solution of which only integrate the plane waves with some unknown directions. 
We find that PWNN works much better than PWPUM at this case. 
Unlike using the plane wave basis with fixed directions in PWPUM, PWNN can learn a group of optimized plane wave basis which can better predict the unknown directions of the solution. 
The proposed approach may provide some new insights in the aspect of applying deep learning in Helmholtz equation.
\end{abstract}


\maketitle

\section{Introduction}
We propose and study a plane wave activation based neural network (PWNN) for solving the Helmholtz boundary value problem
\beq\label{Helm-BV}
\Delta u+k^2u &= 0\quad {\rm in}\ \Omega \\
\frac{\partial u}{\partial n}+ \mathbf{i} ku &= g\quad {\rm on}\ \partial\Omega
\eeq
where $k>0$ is the wavenumber, $g\in L^2(\partial\Omega)$ is the boundary condition and $\Omega \subset \R^2$ is a lipschitz domain. 
We remark that the proposed approach is a basic technique for Helmholtz equation. 
To clarify the discussion, we only consider a simple case of Helmholtz equation, i.e. the Helmholtz boundary value problem \eqref{Helm-BV}. 
We will develop PWNN for some general Helmholtz equations, e.g. unbounded Helmholtz equation, Helmholtz scattering equation, at the future. 

As a famous and significant partial differential equation (PDE), Helmholtz equation \eqref{Helm-BV} appears in diverse scientific and engineering applications, since its solutions can represent the phenomenon of wave propagation, e.g. acoustic wave,  electromagnetic wave, and seismic wave. 
The numerical methods for solving Helmholtz equation have attracted great interests of computational mathematicians, many numerical method like finite element method \cite{Preston1984Finite,ihlenburg1995finite}, finite difference method \cite{newman1995frequency-domain}, boundary element method \cite{Partridge1992The,shen2007an}, plane wave method \cite{babuska1997the,buffa2008error,gittelson2009plane,hiptmair2011plane}, are developed for solving Helmholtz equation in different ways.  

Since deep learning has made great success in engineering, researchers have become increasingly interested in using neural network to solve PDEs. Many relevant papers have been published around different topics in recent years. 
A deep Ritz method is presented by \cite{weinan2017deep}, it can use the variational formulation of elliptic PDEs and use the resulting energy as an objective function for the optimization of the parameters of the neural network.
\cite{zang2019weak} convert the problem of finding the weak solution of PDEs into an operator norm minimization problem induced from
the weak formulation, then parameterized the weak solution and the test function in the weak formulation as the primal and adversarial networks respectively, this method is called weak adversarial network (WAN). 
The physics-informed neural networks (PINN) is presented by \cite{raissi2019physics}, which take initial conditions and boundary conditions as the penalties of the optimization objective loss function.
To deal with the noisy data, \cite{yang2020b-pinns:} proposes the Bayesian Physics-Informed Neural Networks (B-PINN), it use the Hamiltonian
Monte Carlo or the variational inference as a posteriori estimator. Compared with PINNs, B-PINNs obtain more accurate predictions in scenarios
with large noise. 
Unlike using variational formulation, this network needs high-order derivative information, it calculated the high-order differential in the framework of TensorFlow by using the backpropagation algorithm and automatic differentiation. 
The network can give a high accuracy solution for both forward and inverse problems.
\cite{sirignano2018dgm} proposes a deep Galerkin method. Due to the excessive computational cost of high-order differentiation of the neural network, the network bypasses this issue through the Monte Carlo method.
\cite{2020Implicit} proposes sinusoidal representation networks (SIREN), in which sin(x) is used as activation function, and it performs better than other activation function like Sigmoid, Tanh or ReLU in Poisson equation and the Helmholtz and wave equations.
In addition, deep learning tools are used to solve various practical problems related to PDE recently \cite{rudd2013constrained,greydanus2019hamiltonian,berg2018unified}. 
In conclusion, these existing works all focused on the design of the objective of using neural network to solve the general PDEs, not to develop new neural network for a specific PDE.
 
Motivated by the traditional activation based neural network (TANN), with Sigmoid, Tanh or ReLU as acitvation function \cite{raissi2019physics}, we introduce the PWNN method which use the plane wave function $e^{\mathbf{i} {x}}$ as activation function, since the plane wave is the basic component of the solution of Helmholtz equation.
In addition, we find that the PWNN with one hidden layer is a more generalized form of the plane wave partition of unity method ~\cite{babuska1997the,buffa2008error,gittelson2009plane,hiptmair2011plane} (PWPUM). 
In the plane wave method, the direction of the basis function is closely related to the accuracy of the solution.
In PWPUM, the direction of basis function can only be selected averagely.
\cite{2004Numerical,2011Numerical} proposes NMLA to calculate the local direction of the true solution.
A ray-FEM method is presented by \cite{2017Learning}, it first uses the standard finite element method to solve the low-frequency problem, then uses NMLA to calculate the local direction of the low-frequency numerical solution, and finally uses it as the direction of plane wave basis function to solve the high-frequency problem.
\cite{monk1999a} investigate the use of energy function to approximate the Helmholtz equation (LSM), the bases used in the discrete method are either plane waves of Bessel functions.  
\cite{2014A} proposes a wave tracking strategy (LSM-WT), in each element, a same rotation angle is added to all plane wave basis functions with uniformly selected direction, by solving the double minimization problem of energy function in \cite{monk1999a}, the error can be effectively reduced.
Based on this idea, we use neural network to construct this double minimization problem, using Monte Carlo method instead of integration. We call this method PWPUM-WT.   
Unlike using the plane wave basis with fixed directions in PWPUM, 
PWNN may combine several plane waves with learned amplitudes and directions by using plane wave acivation. 
As the direction of each basis function can be rotated by different angles, PWNN is more flexible than PWPUM-WT.
Our experiment show the proposed PWNN method is more efficient than TANN, SIREN and higher accuracy than PWPUM, PWPUM-WT method. 

\section{From TANN to PWNN}
This section introduces plane wave activation based neural network (PWNN) from the view of traditional activation based neural network (TANN).
We first introduce a well known neural network based framework for solving PDEs, i.e. the physics-informed neural networks~\cite{raissi2019physics}. 
The entire neural network consists of $L+1$ layers, 
where layer $0$ is the input layer and layer $L$ is the output layer.
Layers $0<l<L$ are the hidden layers.
All of the layers have an activation function, excluding the output layer.
As we known, the activation function can take the form of Sigmoid, Tanh (hyperbolic tangents) or ReLU (rectified linear units).

We denote $d_0, d_1, ..., d_L$ as a list of integers, with $d_0, d_L$ representing the lengths of input signal and output signal of the neural network. 
Define function ${\bf T}_l: \mathbb{R}^{d_l} \rightarrow \mathbb{R}^{d_{l+1}}$, $0 \leq l < L$,
\beqn
{\bf T}_{l}({\bf x})={\bf W}_{l} {\bf x} + {\bf b}_{l},
\eeqn
where ${\bf W}_{l} \in \mathbb{R}^{d_{l+1} \times d_{l}}$ and ${\bf b}_{l} \in \mathbb{R}^{d_{l+1}}$. 
Thus, we can simply represent a deep fully connected feedforward neural network using the composite function $h(\cdot;\theta):\mathbb{R}^{d_0} \rightarrow \mathbb{R}^{d_L}$,
\beqn
h(\cdot;\theta) = {\bf T}_{L-1}\circ\sigma\circ {\bf T}_{L-2} \circ \cdots \circ {\bf T}_1 \circ \sigma \circ {\bf T}_0,
\eeqn
where $\sigma$ is the activation function and $\theta:=\{{\bf W}_{l},{\bf b}_{l}\}$ represents the collection of all parameters. 
We only consider $d_0 = 2, d_L = 1$, a two-dimensional scalar problem for Helmholtz equation \eqref{Helm-BV}. 

Solving a general PDE defined as
\beq \label{eq:PDE}	
\mathcal{L} (u) &= f, \quad \text{in}\, \Omega, \\
\mathcal{B} (u) &= g, \quad \text{on}\, \partial\Omega, 
\eeq 
by DNN is a physics-informed minimization problem with the objective $\mathcal{M}(\theta)$ consisting of two terms as follows:
\beq\label{eq:objective}
\theta^*={\rm argmin}_{\theta}\, \mathcal{M}(\theta)&:=\mathcal{M}_{\Omega}(\theta)+\lambda\mathcal{M}_{\partial \Omega}(\theta), \\
\mathcal{M}_{\Omega}(\theta):=\frac{1}{N_f}\sum_{i=1}^{N_f}
\left\vert \mathcal{L}\left(h({\bf x}_f^{i};\theta)\right)-f({\bf x}_f^i)\right\vert^{2},
&\mathcal{M}_{\partial \Omega}(\theta):=\frac{1}{N_g}\sum_{i=1}^{N_g} 
\left\vert \mathcal{B}\left(h({\bf x}_g^{i};\theta) \right)-g({\bf x}_g^{i}) \right\vert^{2}, 
\eeq
where $\{{\bf x}_f^i \}_{i=1}^{N_f}$ and $\{{\bf x}_g^i \}_{i=1}^{N_g}$ are the collocation points in the inside and on the boundary of domain $\Omega$, respectively. $\lambda$ is a hyperparameter, by adjusting $\lambda$ to make the internal and boundary initial errors of the same order, usually there will be a better result.
The domain term $\mathcal{M}_{\Omega}$ and boundary term $\mathcal{M}_{\partial \Omega}$ enforce the condition that the desired optimized neural network $h(\cdot;\theta^*)$ satisfies $\mathcal{L}(u)=f$ and $\mathcal{B}(u)=g$, respectively. 
The effective and efficient stochastic gradient descents~\cite{lecun2012efficient} with minibatches are recommended to solve the optimization problem \eqref{eq:objective}.


For Helmholtz equation, since the solution is with complex value, we make a slight modification on the network architecture by redefining ${\bf W}_{L-1} \in \mathbb{C}^{d_{L} \times d_{L-1}},{\bf b}_{L-1} \in \mathbb{C}^{d_{L}}$. 
We may consider directly use the complex neural network~\cite{trabelsi2018deep} at the future. 
By this definition, we have $h(\cdot;\theta):\mathbb{R}^{d_0} \rightarrow \mathbb{C}^{d_L}$. 
Now, we let $\mathcal{L}=\Delta+k^2, \mathcal{B}=\partial/\partial n+ \mathbf{i}k$ and $f=0$, 
and let the activation function $\sigma(x)$ be ${\rm tanh}(x)$ or $e^{\mathbf{i}x}$ in minimization problem \eqref{eq:objective} to produce a \textbf{TANN} or \textbf{PWNN} based numerical method for solving Helmholtz equation \eqref{Helm-BV} respectively. 

As the key issue of numerically solving PDEs, a basic and interesting question is as follows: is it possible to use $h({\bf x};\theta)$ to approximate the solution of PDE? 
A well-known answer is the universal approximation property (UAP): if the solution $u({\bf x})$ is bounded and continuous, then $h({\bf x};\theta)$ can approximate $u({\bf x})$ to any desired accuracy, given the increasing hidden neurons~\cite{hornik1989multilayer, cybenko1989approximation, hornik1991approximation}. 
We should notice that this conclusion required a simple assumption on the activation function, i.e. $\sigma(x)$ should be bounded and monotonically-increasing function. 
However, it's easy to check that the plane wave activation doesn't satisfy this assumption. 
We open a theoretical problem that whether PWNN has the UAP at the future research. 

\section{From PWNN to PWPUM}
In this section , we will try to bridges PWNN to the plane wave partition of unity method (PWPUM), which is the famous and classical numerical method with complexity and approximation analysis. 
If we consider a simpler PWNN with one hidden layer, let $\mathbf{T}_0:=(w_1,w_2,\cdots,w_{d_1})^T$, 
$\mathbf{b}_0:=(b_1,b_2,\cdots,b_{d_1})^T$, $\mathbf{T}_1:=(s_1,s_2,\cdots,s_{d_1})$, 
and $w_i \in \R^{2}$, $b_i \in \R$, $s_i \in \mathbb{C}$, 
we will have the output function with the form
\beq\label{PWNN-H1}
h(\mathbf{x};\theta) = \sum_{i=1}^{d_1} s_i e^{\mathbf{i} (w_i^T \mathbf{x}+b_i)}. 
\eeq
Notice that $e^{\mathbf{i}b_i}$ can be equivalently learned by learning $s_i$, 
we can ignore $b_i$ in \eqref{PWNN-H1} to obtain an equivalent output formulation
\beq\label{PWNN-H1-v2}
h(\mathbf{x};\theta) = \sum_{i=1}^{d_1} s_i e^{\mathbf{i} w_i^T \mathbf{x}},
\eeq
with $\theta:=\{s_i,w_i\}_{i=1}^{d_1}$. So $\theta$ of PWNN with one hidden layer is determined by solving the following problem,
\beq\label{op_PWNN}
{\rm min}_{\theta}\, \int_{\Omega}|\mathcal{L}(h(\mathbf{x};\theta))-f(\mathbf{x})|^2 + \lambda\int_{\partial\Omega}|\mathcal{B}(h(\mathbf{x};\theta))-g(\mathbf{x})|^2
\eeq 
Here we write the loss function in integral form. In fact, we use Monte Carlo method instead of integration in neural network. In the standard finite element method or other methods, there will also be errors in numerical integration. This part of the error is only a small part when not pursuing too high accuracy, e.g. relative $L^2$ error less than 1e-8.

Now we come to introduce the plane wave partition of unity method ~\cite{babuska1997the,buffa2008error,gittelson2009plane,hiptmair2011plane} (\textbf{PWPUM}). 
For the Helmholtz boundary value equation with wavenumber $k$, $1\leq i\leq d_1$, 
let 
\beq\label{def:k_i}
k_i:=(k\cos \frac{2\pi i}{d_1},k\sin \frac{2\pi i}{d_1}), 
\eeq
the PWPUM is to find a numerical solution 
with the form
\beq\label{PWPUM}
u(\mathbf{x};\beta) = \sum_{i=1}^{d_1} s_i e^{\mathbf{i} k_i^T \mathbf{x}},
\eeq
where we define $\beta:=\{s_i\}_{i=1}^{d_1}$ the parameters which need to be optimized. 
We define the discrete function space with form \eqref{PWPUM} depending on $\beta$ as $V_h$. 
Instead of directly solving an physics-informed optimization problem \eqref{eq:objective} in PWNN, 
PWPUM instead find $u\in V_h$, such that, $\forall\, v \in V_h$,
\beq\label{PWPUM-V}
\int_{\Omega}\,(\Delta u+k^2u)\overline{v}\,dV &=0, \\
\int_{\partial\Omega}\,(\frac{\partial u}{\partial n}+\mathbf{i}ku)\overline{v}\,dS&=\int_{\partial\Omega}\,g\overline{v}\,dS.
\eeq
By the formulation of $u$ in \eqref{PWPUM} and the definition of $k_i$ in \eqref{def:k_i}, 
we have $\Delta u+k^2u=0$, yielding that the first volume integral term of \eqref{PWPUM-V} vanishes. 
Since the Helmholtz operator and the discrete function space are both linear, we can  
take the test function $v=e^{\mathbf{i} k_i^T \mathbf{x}}$, $1 \leq i \leq d_1$. 
Then, to find $u\in V_h$ satisfying \eqref{PWPUM-V} is to solve a linear system $\mathbf{M} \beta = \mathbf{G}$ with $\mathbf{M} \in \mathbb{C}^{d_1\times d_1}, \mathbf{G} \in \mathbb{C}^{d_1}$,and
\beqn
\mathbf{M}_{i,j} = \int_{\partial\Omega}(\frac{\partial e^{\mathbf{i} k_i^T \mathbf{x}}}{\partial n}+\mathbf{i}ke^{\mathbf{i} k_i^T \mathbf{x}})&{e^{-\mathbf{i} k_j^T \mathbf{x}}}dS, \quad
\mathbf{G}_i = \int_{\partial\Omega}ge^{-\mathbf{i} k_i^T \mathbf{x}}dS. 
\eeqn

Actually, one may check that the above linear system is equivalent to solving the following problem:
\beq\label{op_PWPUM}
{\rm min}_{\beta}\, \int_{\Omega}|\mathcal{L}(u(\mathbf{x};\beta))-f(\mathbf{x})|^2 + \lambda\int_{\partial\Omega}|\mathcal{B}(u(\mathbf{x};\beta))-g(\mathbf{x})|^2.
\eeq
In fact, the basis function of PWPUM is plane wave function and the wave number is strictly equal to $k$. The internal integral error is strictly equal to 0, so this optimization problem is independent of $\lambda$.

For the PWPUM method, the direction of the plane wave basis function is uniformly selected because the prior information of the solution is not known. Because any direction may be the direction of the exact solution, PWPUM may need a lot of basis functions to achieve the required accuracy.
In order for PWPUM to use the direction that can be changed, \cite{2014A} proposes a wave tracking strategy (LSM-WT), in each element, a same rotation angle is added to all plane wave basis functions with uniformly selected direction. 
Here we can also use this wave tracking strategy in one element, the form of solution \eqref{PWPUM} becomes:
\beq\label{PWPUM-WT}
u(\mathbf{x};\beta,\alpha) = \sum_{i=1}^{d_1} s_i e^{\mathbf{i} k_i^T(\alpha) \mathbf{x}},
\eeq
where
\beqn
\beta:=\{s_i\}_{i=1}^{d_1},\ k_i^T(\alpha)=\begin{bmatrix} \cos\alpha & \sin\alpha \\ -\sin\alpha & \cos\alpha \end{bmatrix}k_i^T,
\eeqn
and the optimization problem \eqref{op_PWPUM} becomes:
\beq\label{op_PWPUM-WT}
{\rm min}_{\beta,\alpha}\, \int_{\Omega}|\mathcal{L}(u(\mathbf{x};\beta,\alpha))-f(\mathbf{x})|^2 + \lambda\int_{\partial\Omega}|\mathcal{B}(u(\mathbf{x};\beta,\alpha))-g(\mathbf{x})|^2.
\eeq

In addition, if we limit $w_i$ in \ref{PWNN-H1-v2} to $w_i^T=k_i^T(\alpha)$, 
then this is the form of PWPUM-WT, the output formulation $h(\mathbf{x};\theta)$ is equal to $u(\mathbf{x};\beta,\alpha)$ in \eqref{PWPUM-WT},
with $\theta:=\{s_i,\alpha\}_{i=1}^{d_1}$.
 
Apparently, comparing PWNN \eqref{PWNN-H1-v2}, PWPUM \eqref{PWPUM} and PWPUM-WT \eqref{PWPUM-WT} will show that PWPUM-WT is a generalization of PWPUM. Furthermore, PWNN is a generalization of PWPUM-WT.
Unlike using the plane wave basis with fixed directions \eqref{def:k_i} in PWPUM, 
PWNN will learn a group of optimized plane waves $e^{\mathbf{i} w_i^T \mathbf{x}}$ in \eqref{PWNN-H1-v2} with learned both amplitude $||w_i||$ and direction ${\rm arccos}(w_{i,1}/||w_i||)$.
Since the approximation analysis of PWPUM was well addressed, see Theorem 3.9 of \cite{hiptmair2011plane}, and the optimization problem of PWNN \eqref{op_PWNN} is obviously a generalization of that in PWPUM \eqref{op_PWPUM}, it's trivial to prove that the approximation error of PWNN with one hidden layer is no more than that of PWPUM theoretically.
Actually, we can only consider learning the directions by adding a regularization term, like $\sum_{i=1}^{d_1}(||w_i||-k)^2$ in the minimization objective \eqref{eq:objective}. 
However, motivated by the dispersion correction technique~\cite{cocquet2017a,cocquet2019a}, i.e. using different wavenumber in the discretization scheme may reduce the pollution effect~\cite{deraemaeker1999dispersion} of the Helmholtz problems with high wavenumber, we ignore adding this regularization term.

\section{Experiments}
We experimental investigate the performance of PWNN from two aspects: the exact solution of equation \eqref{Helm-BV} is a combination of plane waves with all known directions (KD) or unknown directions (UD). 
Apparently, the UD problem is more common in practical problems. 
For instance, if we consider a Helmholtz scattering equation with complicated scatters, i.e. the propagation of acoustic wave in a complicated situation, we would not expect to accurately predict the directions of the multi reflection waves in advance.
Considering that the sin activation based neural network (SIREN) performs better than the traditional activation based neural network (TANN) in some cases, we add the comparison between PWNN and SIREN in the experiments. 
The experimental results on problems KD and UD demonstrate that: PWNN works much better than TANN and SIREN, i.e. the introduction of plane wave as activation unit is indeed helpful for Helmholtz equation; 
PWPUM-WT has limited advantages over PWPUM, i.e. adding only one degree of freedom of rotation angle is not enough.
PWNN has competitive performance with PWPUM under problem KD but outperforms PWPUM under problem UD, i.e. PWNN is indeed a generalization of PWPUM but more practical.
We introduce the experimental setting first. 

\textbf{Experimental setting}. 
We simply set the network architectures as having equal units for each layer. 
The $Layers$ used later indicates the hidden layers in the network, i.e. $Layers:=L-1$. 
The $Units$ indicates the units per layer, i.e. $d_1=d_2=\cdots=d_{L-1}=Units$. We chose Tanh as activation function in TANN. In both TANN and PWNN, we use limited-memory BFGS algorithm \cite{nocedal1980updating} to update parameters. L-BFGS is an improved algorithm for quasi-Newton method: BFGS. In BFGS, the  approximate Hesse matrix $B_k^{-1}$ is stored at every step, which wastes a lot of storage space in high dimensional cases. In L-BFGS, only the recent $m$ steps' iterative information is saved for calculation $B_k^{-1}$ to reduce the storage space of data. Here we set $m=50$. The stop criterion of the inner iteration is $\|\nabla\mathcal{M}\|_{\infty}<2 \times 10^{-16}$, or that we exceeded the maximum number of allowed iterations, set as 50000 here. $N_f,N_g$ denotes the number of training points in the inside and on the boundary. The training points in the inside are randomly selected, while the training points on the boundary are uniformly selected. 
The test points used to estimate the relative $L^2$ errors are uniformly sampled by row and column, with 10000 in all. 
We remark that this paper ignores the detailed discussion of computational cost, 
since all the numerical examples only take up to several minutes to get satisfying results. 
We may consider comparing the computational cost of TANN, PWNN, and PWPUM for a 3D Helmholtz problem at the future. 

\subsection{Problem with known directions (KD)}

We consider a square domain $\Omega=[0,1]\times[-0.5,0.5]$. 
Let the analytical solutions are the circular waves given, in polar coordinates $\mathbf{x}=(r\cos\theta ,r\sin\theta)$, by
\beq\label{es1}
u_*(\mathbf{x})=J_{\xi}(kr)e^{i\xi\theta},
\eeq
where $J_{\xi}$ denotes the Bessel function of the first kind and order\ $\xi$, we choose $\xi=1$.
It's easy to check that $\Delta u_* + k^2u_* = 0$. 
We substitute $u_*$ into \eqref{Helm-BV} to compute $g$.
We denote the numerical solution by $u_h$ and define:
\beqn
\varepsilon :=\frac{\|(|u_*|-|u_h|)\|_2}{\|(|u_*|)\|_2},
\eeqn
which denote relative $L^2$ error between predicted and the exact solution in length, then we define the \textbf{accuracy} as:
\beqn
Accuracy := -\lg\varepsilon.
\eeqn

We present the results of solving Helmholtz equation \eqref{Helm-BV} by TANN, SIREN and PWNN. 
Table \ref{TANN_arch,k=5},\ref{TANN_arch,k=10},\ref{TANN_arch,k=20} shows the relative $L^2$ error between predicted and the exact solution for different network architectures, while the total numbers of training points are fixed to $N_g=5k\times 4$ and $N_f=5k^2$, $k$ denotes the wave number of exact solution. 
Here we take $Layers=1,2,3,4$ and $Units=k,2k,4k$ respectively for all these networks. We choose $Units,\ N_g$ and $N_f$ increase with $k$ in order to improve the approximation of the model and reduce the generalization error to deal with the problem of higher $k$.
As expected, when the numbers of $Layers$ is fixed, more $Units$ leads to smaller relative $L^2$ errors on the whole. Also, for a fixed number of $Units$, relative errors is likely to decreases with the increase of $Layers$. The accuracy acquired by PWNN is higher than that of SIREN and TANN with the same numbers of $Layers$ and $Units$, especially for $Layers=1$. With the increase of wave number $k$, the accuracy of PWNN is improved more. When $Units<2k$, a multi-layer PWNN achieves higher accuracy than one hidden layer. But for larger $Units$, one hidden layer performs better. This means that the approximation ability of one hidden layer PWNN is enough at this case, increasing the number of $Layers$ will make network parameters difficult to optimize.

\begin{table}[htbp]
	\centering
	\caption{\textbf{KD problem}. Relative $L^2$ error $\varepsilon$ (Accuracy improvement compared with TANN) between predicted and the exact solution for different numbers of hidden layers and different numbers of units per layer in TANN, SIREN and PWNN. Here, the total numbers of training data are fixed to $N_g=25\times 4$ and $N_f=125$. The wave number in exact solution \ref{es1} is taken as $k=5$.}
	\label{TANN_arch,k=5}
	\small{
	\begin{tabular}{c|c|ccc}
		\toprule
		network & \diagbox{$Layers$}{$Units$} & 5 & 10 & 20 \\
		\midrule
		\multirow{4}*{TANN}& 1 & 1.8e-1 & 8.2e-3 & 4.5e-3 \\
		
		&2 & 9.5e-2 & 3.4e-3 & 1.8e-3 \\
		
		&3 & 2.0e-2 & 2.3e-3 & 2.0e-3 \\
		
		&4 & 1.6e-2 & 2.8e-3 & 3.0e-3 \\
		\midrule
		\multirow{4}*{SIREN}& 1 & 8.5e-2($+$0.33) & 6.5e-3($+$0.10) & 8.3e-3($-$0.27) \\
		
		&2 & 1.4e-2($+$0.83) & 1.7e-3($+$0.30) & 3.2e-3($-$0.25) \\
		
		&3 & 6.9e-3($+$0.46) & 3.0e-3($-$0.12) & 1.7e-3($+$0.07) \\
		
		&4 & 5.6e-3($+$0.46) & 2.0e-3($+$0.15) & 1.1e-3($+$0.44) \\
		\midrule
		\multirow{4}*{PWNN}& 1 & 1.8e-1($+$0.00) & 4.1e-4($+$1.30) & 6.9e-6($+$2.81) \\
		
		&2 & 1.8e-3($+$1.72) & 4.7e-4($+$0.86) & 4.4e-4($+$0.61) \\
		
		&3 & 2.8e-3($+$0.85) & 2.9e-4($+$0.90) & 5.9e-4($+$0.53) \\
		
		&4 & 1.8e-3($+$0.95) & 4.1e-4($+$0.83) & 6.7e-4($+$0.65) \\
		\bottomrule
	\end{tabular}
	}
\end{table}

\begin{table}[htbp]
	\centering
	\caption{\textbf{KD problem}. Relative $L^2$ error $\varepsilon$ (Accuracy improvement compared with TANN) between predicted and the exact solution for different numbers of hidden layers and different numbers of units per layer in TANN, SIREN and PWNN. Here, the total numbers of training data are fixed to $N_g=50\times 4$ and $N_f=500$. The wave number in exact solution \ref{es1} is taken as $k=10$.}
	\label{TANN_arch,k=10}
	\small{
	\begin{tabular}{c|c|ccc}
		\toprule
		network & \diagbox{$Layers$}{$Units$} & 10 & 20 & 40 \\
		\midrule
		\multirow{4}*{TANN}& 1 & 4.3e-1 & 8.2e-2 & 2.6e-2 \\
		
		&2 & 3.9e-2 & 7.1e-3 & 4.0e-3 \\
		
		&3 & 2.4e-2 & 4.7e-3 & 3.9e-3 \\
		
		&4 & 2.1e-2 & 4.3e-3 & 2.9e-3 \\
		\midrule
		\multirow{4}*{SIREN}& 1 & 1.0e-1($+$0.63) & 5.7e-2($+$0.16) & 3.1e-2($-$0.08) \\
		
		&2 & 5.7e-3($+$0.84) & 4.7e-3($+$0.09) & 4.1e-3($-$0.01) \\
		
		&3 & 3.1e-3($+$0.89) & 2.6e-3($+$0.26) & 3.4e-3($+$0.06) \\
		
		&4 & 5.0e-3($+$0.62) & 2.5e-3($+$0.24) & 1.2e-3($+$0.38) \\
		\midrule
		\multirow{4}*{PWNN}& 1 & 1.5e-1($+$0.46) & 1.7e-5($+$3.68) & 5.4e-6($+$3.68) \\
		
		&2 & 1.1e-3($+$1.55) & 7.3e-4($+$0.99) & 4.0e-4($+$1.00) \\
		
		&3 & 9.3e-4($+$1.41) & 2.4e-4($+$1.29) & 1.2e-3($+$0.51) \\
		
		&4 & 4.2e-3($+$0.70) & 1.4e-4($+$1.49) & 9.1e-4($+$0.50) \\
		\bottomrule
	\end{tabular}
	}
\end{table} 

\begin{table}[htbp]
	\centering
	\caption{\textbf{KD problem}. Relative $L^2$ error $\varepsilon$ (Accuracy improvement compared with TANN) between predicted and the exact solution for different numbers of hidden layers and different numbers of units per layer in TANN, SIREN and PWNN. Here, the total numbers of training data are fixed to $N_g=100\times 4$ and $N_f=2000$. The wave number in exact solution \ref{es1} is taken as $k=20$.}
	\label{TANN_arch,k=20}
	\small{
	\begin{tabular}{c|c|ccc}
		\toprule
		network & \diagbox{$Layers$}{$Units$} & 20 & 40 & 80 \\
		\midrule
		\multirow{4}*{TANN}& 1 & 6.0e-1 & 3.2e-1 & 1.4e-1 \\
		
		&2 & 6.8e-2 & 3.8e-2 & 3.3e-2 \\
		
		&3 & 3.5e-2 & 9.5e-3 & 1.0e-2 \\
		
		&4 & 1.9e-2 & 1.4e-2 & 9.2e-3 \\
		\midrule
		\multirow{4}*{SIREN}& 1 & 5.2e-1($+$0.06) & 2.5e-2($+$1.10) & 2.8e-1($-$0.30) \\
		
		&2 & 2.4e-1($-$0.55) & 1.2e-2($+$0.50) & 2.5e-2($+$0.12) \\
		
		&3 & 1.2e-2($+$0.46) & 9.4e-3($+$0.00) & 1.5e-2($-$0.18) \\
		
		&4 & 9.6e-3($+$0.30) & 1.2e-2($+$0.07) & 5.0e-3($+$0.26) \\
		\midrule
		\multirow{4}*{PWNN}& 1 & 1.8e-1($+$0.52) & 7.9e-7($+$5.60) & 4.7e-6($+$4.47) \\
		
		&2 & 1.4e-2($+$0.69) & 1.2e-3($+$1.50) & 1.0e-3($+$1.52) \\
		
		&3 & 1.7e-3($+$1.31) & 5.2e-4($+$1.26) & 1.5e-3($+$0.82) \\
		
		&4 & 1.7e-3($+$1.04) & 2.2e-4($+$1.80) & 1.0e-3($+$0.96) \\
		\bottomrule
	\end{tabular}
	}
\end{table}

\begin{table}[htbp]
	\centering
	\caption{\textbf{KD problem}. Average ($\pm$Standard Deviation) of the relative $L^2$ error between predicted and the exact solution for different numbers of training data in TANN, SIREN and PWNN. Here we randomly generate 50 groups of internal sample points, calculate the average value and standard deviation of relative error. The network architecture of TANN, SIREN and PWNN is fixed to $Layers=1,Units=40$. The wave number in exact solution \ref{es1} is taken as $k=20$.}
	\label{TANN_data1}
	\small{
	\resizebox{\textwidth}{!}{
	\begin{tabular}{c|c|cccc}
		\toprule
		network & \diagbox{$N_g$}{$N_f$} &200 & 500 & 1000 & 2000 \\
		\midrule
		\multirow{4}*{TANN}& 20$\times$4 & 4.5e-1($\pm$5.1e-2) & 4.2e-1($\pm$7.5e-2) & 3.6e-1($\pm$8.0e-2) & 3.9e-1($\pm$6.7e-2) \\
		
		&50$\times$4 & 4.5e-1($\pm$5.4e-2) & 4.0e-1($\pm$6.7e-2) & 3.8e-1($\pm$8.6e-2) & 3.2e-1($\pm$8.6e-2) \\
		
		&100$\times$4 & 4.4e-1($\pm$6.1e-2) & 3.9e-1($\pm$6.6e-2) & 3.4e-1($\pm$7.4e-2) & 3.2e-1($\pm$8.3e-2) \\
		
		&200$\times$4 & 4.5e-1($\pm$5.7e-2) & 3.8e-1($\pm$8.1e-2) & 3.4e-1($\pm$9.0e-2) & 3.1e-1($\pm$9.4e-2) \\
		\midrule
		\multirow{4}*{SIREN}& 20$\times$4 & 5.3e-1($\pm$8.1e-2) & 5.2e-1($\pm$9.5e-2) & 5.2e-1($\pm$1.3e-1) & 4.8e-1($\pm$1.3e-1) \\
		
		&50$\times$4 & 5.0e-1($\pm$1.0e-1) & 4.8e-1($\pm$1.3e-1) & 4.7e-1($\pm$1.4e-1) & 4.8e-1($\pm$1.4e-1) \\
		
		&100$\times$4 & 4.9e-1($\pm$1.3e-1) & 4.6e-1($\pm$1.4e-1) & 4.7e-1($\pm$1.4e-1) & 4.3e-1($\pm$1.5e-1) \\
		
		&200$\times$4 & 5.1e-1($\pm$6.8e-2) & 5.0e-1($\pm$9.2e-2) & 4.9e-1($\pm$1.0e-1) & 4.0e-1($\pm$1.7e-1) \\
		\midrule
		\multirow{4}*{PWNN}& 20$\times$4 & 7.7e-7($\pm$2.2e-7) & 7.4e-7($\pm$2.5e-7) & 6.9e-7($\pm$1.6e-7) & 6.8e-7($\pm$1.8e-7) \\
		
		&50$\times$4 & 7.4e-7($\pm$2.1e-7) & 6.6e-7($\pm$2.1e-7) & 6.6e-7($\pm$1.8e-7) & 6.5e-7($\pm$2.2e-7) \\
		
		&100$\times$4 & 7.9e-7($\pm$2.6e-7) & 6.7e-7($\pm$1.9e-7) & 6.5e-7($\pm$2.2e-7) & 6.3e-7($\pm$1.8e-7) \\
		
		&200$\times$4 & 7.7e-7($\pm$2.6e-7) & 7.0e-7($\pm$1.6e-7) & 6.5e-7($\pm$1.9e-7) & 6.7e-7($\pm$1.9e-7) \\
		\bottomrule
	\end{tabular}}
	}
\end{table}

\begin{table}[htbp]
	\centering
	\caption{\textbf{KD problem}. Average ($\pm$Standard Deviation) of the relative $L^2$ error between predicted and the exact solution for different numbers of training data in TANN, SIREN and PWNN. Here we randomly generate 50 groups of internal sample points, calculate the average value and standard deviation of relative error. The network architecture of TANN, SIREN and PWNN is fixed to $Layers=2,Units=40$. The wave number in exact solution \ref{es1} is taken as $k=20$.}
	\label{TANN_data2}
	\small{
	\resizebox{\textwidth}{!}{
	\begin{tabular}{c|c|cccc}
		\toprule
		network & \diagbox{$N_g$}{$N_f$} &200 & 500 & 1000 & 2000 \\
		\midrule
		\multirow{4}*{TANN}& 20$\times$4 & 1.6e-1($\pm$8.2e-2) & 6.6e-2($\pm$2.2e-2) & 4.6e-2($\pm$1.2e-2) & 4.1e-2($\pm$1.2e-2)\\
		
		&50$\times$4 & 1.5e-1($\pm$7.8e-2) & 5.7e-2($\pm$1.6e-2) & 4.5e-2($\pm$1.2e-2) & 3.5e-2($\pm$9.7e-3) \\
		
		&100$\times$4 & 1.6e-1($\pm$1.0e-1) & 5.7e-2($\pm$1.7e-2) & 3.8e-2($\pm$1.2e-2) & 3.3e-2($\pm$1.0e-2) \\
		
		&200$\times$4 & 1.4e-1($\pm$5.3e-2) & 5.4e-2($\pm$2.0e-2) & 4.2e-2($\pm$1.9e-2) & 3.3e-2($\pm$9.5e-3) \\
		\midrule
		\multirow{4}*{SIREN}& 20$\times$4 & 6.9e-2($\pm$3.7e-2) & 3.7e-2($\pm$1.2e-2) & 2.8e-2($\pm$7.0e-3) & 2.7e-2($\pm$9.0e-3) \\
		
		&50$\times$4 & 6.0e-2($\pm$3.5e-2) & 3.4e-2($\pm$1.0e-2) & 2.3e-2($\pm$6.7e-3) & 2.4e-2($\pm$7.0e-3) \\
		
		&100$\times$4 & 5.7e-2($\pm$2.8e-2) & 3.0e-2($\pm$8.0e-3) & 2.2e-2($\pm$5.7e-3) & 1.9e-2($\pm$6.6e-3) \\
		
		&200$\times$4 & 5.5e-2($\pm$2.1e-2) & 2.8e-2($\pm$8.7e-3) & 2.2e-2($\pm$5.7e-3) & 1.9e-2($\pm$5.8e-3) \\
		\midrule
		\multirow{4}*{PWNN}& 20$\times$4 & 2.1e-1($\pm$2.3e-1) & 2.5e-3($\pm$5.7e-4) & 1.8e-3($\pm$5.0e-4) & 1.4e-3($\pm$7.6e-4) \\
		
		&50$\times$4 & 2.2e-1($\pm$2.4e-1) & 2.5e-3($\pm$5.5e-4) & 1.7e-3($\pm$4.0e-4) & 1.4e-3($\pm$3.3e-4) \\
		
		&100$\times$4 & 2.2e-1($\pm$2.2e-1) & 2.5e-3($\pm$4.9e-4) & 1.8e-3($\pm$4.0e-4) & 1.3e-3(($\pm$2.9e-4)) \\
		
		&200$\times$4 & 2.3e-1($\pm$2.5e-1) & 2.5e-3($\pm$4.7e-4) & 1.8e-3($\pm$3.5e-4) & 1.4e-3($\pm$4.1e-4) \\
		\bottomrule
	\end{tabular}}
	}
\end{table} 

\begin{figure*}[htbp]
	\centering
	\subfigure[$Layers=1,\ k=5$]{
		\begin{minipage}[t]{0.25\textwidth}
			\centering
			\includegraphics[width=1\textwidth]{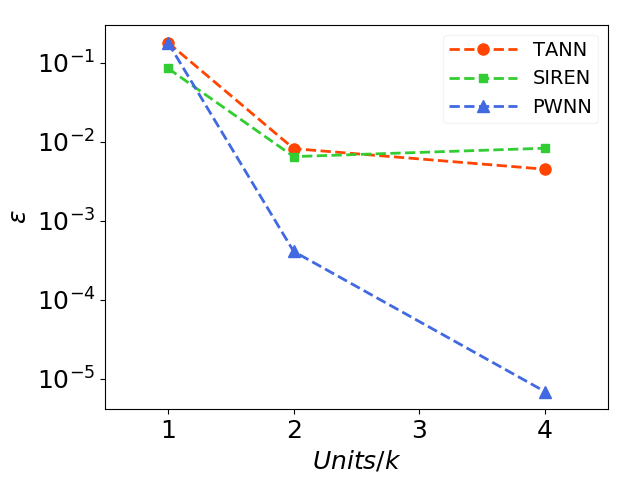}
	\end{minipage}}
	\subfigure[$Layers=1,\ k=10$]{
		\begin{minipage}[t]{0.25\textwidth}
			\centering
			\includegraphics[width=1\textwidth]{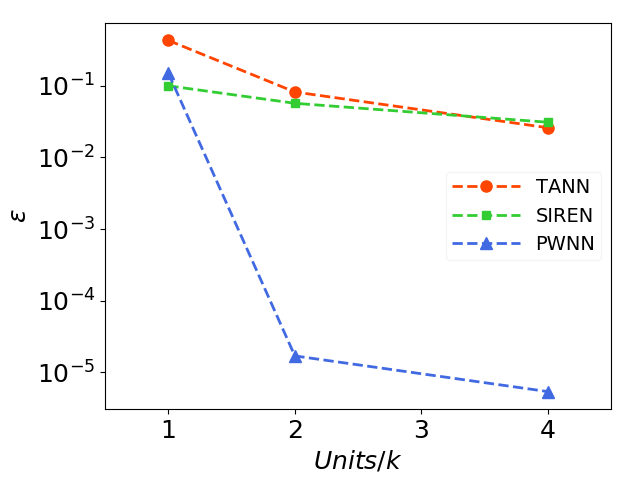}
	\end{minipage}}
    \subfigure[$Layers=1,\ k=20$]{
    	\begin{minipage}[t]{0.25\textwidth}
    		\centering
    		\includegraphics[width=1\textwidth]{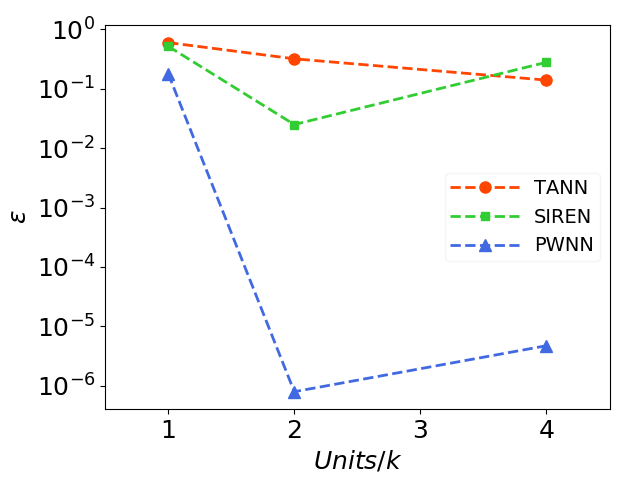}
    \end{minipage}}

    \subfigure[$Layers=2,\ k=5$]{
    	\begin{minipage}[t]{0.25\textwidth}
    		\centering
    		\includegraphics[width=1\textwidth]{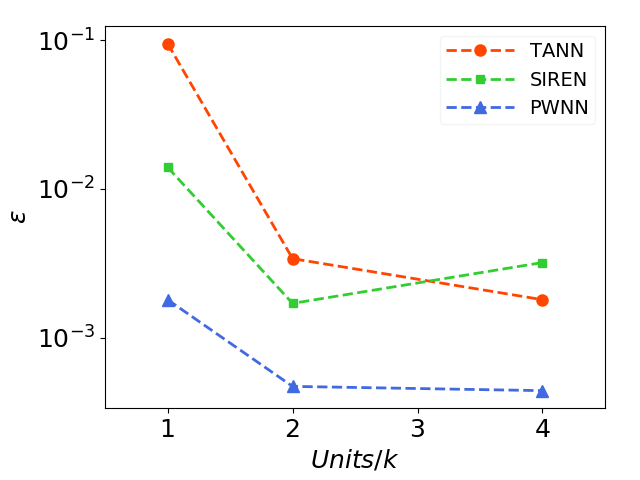}
    \end{minipage}}
    \subfigure[$Layers=2,\ k=10$]{
    	\begin{minipage}[t]{0.25\textwidth}
    		\centering
    		\includegraphics[width=1\textwidth]{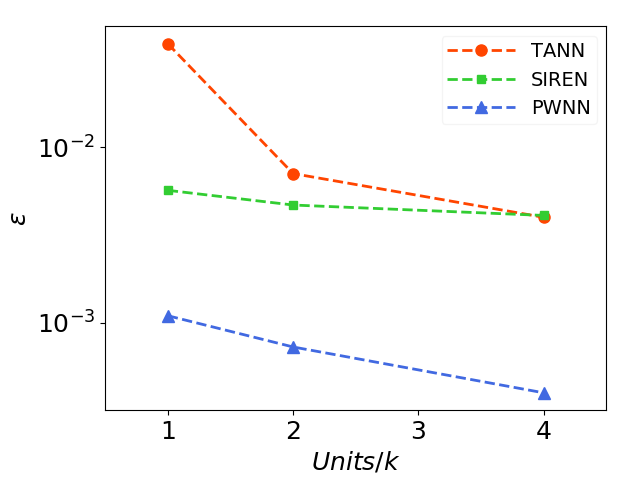}
    \end{minipage}}
    \subfigure[$Layers=2,\ k=20$]{
    	\begin{minipage}[t]{0.25\textwidth}
    		\centering
    		\includegraphics[width=1\textwidth]{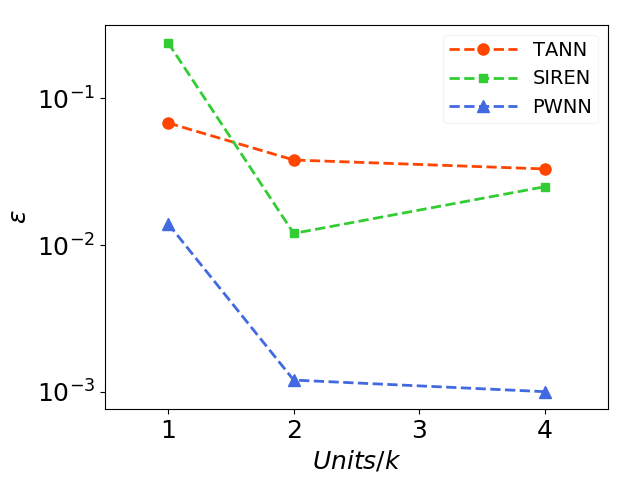}
    \end{minipage}}
	\caption{\textbf{KD problem}. Relative error of TANN, SIREN and PWNN for different $Units$ and wave number $k$, while $k$ is fixed in each subfigure. Here, $Layers$ are respectively taken as 1, 2. The x-axis represents $Units/k$. The total numbers of training points are fixed to $N_g=5k\times 4$ and $N_f=5k^2$.}
	\label{error1}
\end{figure*}

\begin{figure*}[htbp]
	\centering
	\subfigure[$Layers=1$, TANN]{
		\begin{minipage}[t]{0.25\textwidth}
			\centering
			\includegraphics[width=1\textwidth]{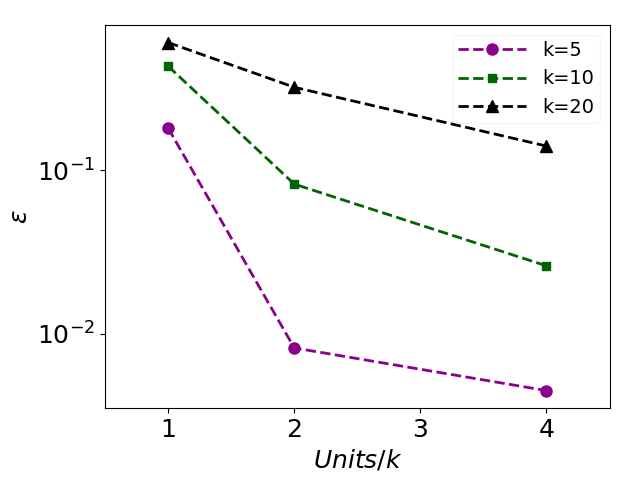}
	\end{minipage}}
	\subfigure[$Layers=1$, SIREN]{
		\begin{minipage}[t]{0.25\textwidth}
			\centering
			\includegraphics[width=1\textwidth]{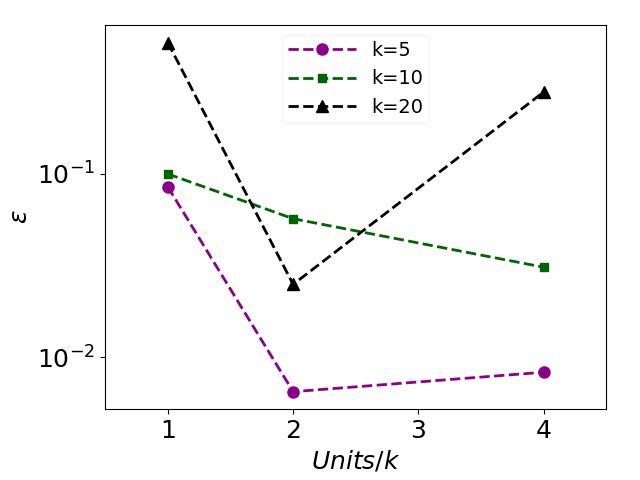}
	\end{minipage}}
	\subfigure[$Layers=1$, PWNN]{
		\begin{minipage}[t]{0.25\textwidth}
			\centering
			\includegraphics[width=1\textwidth]{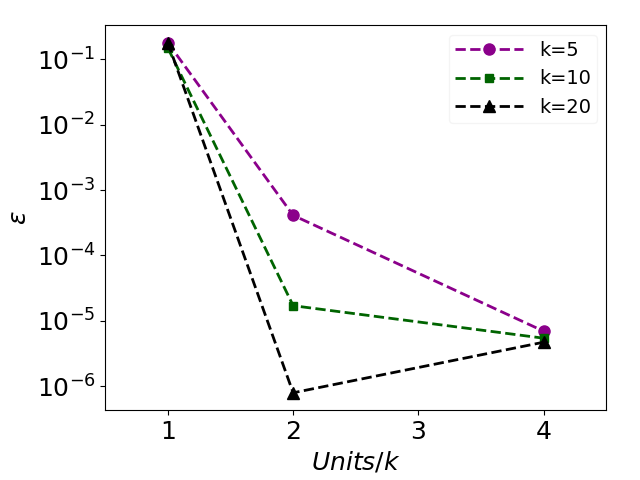}
	\end{minipage}}
	
	\subfigure[$Layers=2$, TANN]{
		\begin{minipage}[t]{0.25\textwidth}
			\centering
			\includegraphics[width=1\textwidth]{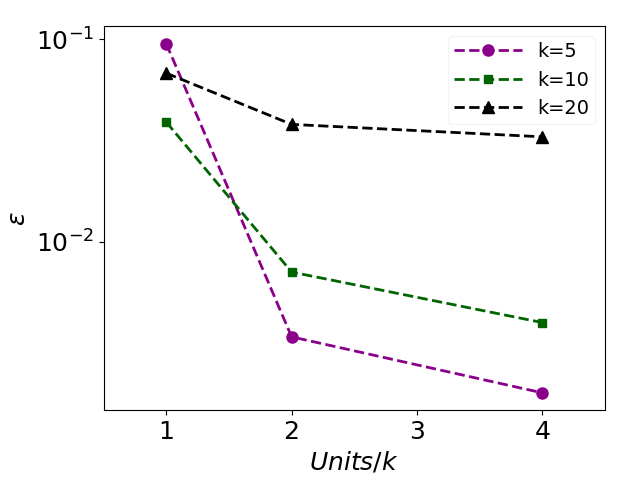}
	\end{minipage}}
	\subfigure[$Layers=2$, SIREN]{
		\begin{minipage}[t]{0.25\textwidth}
			\centering
			\includegraphics[width=1\textwidth]{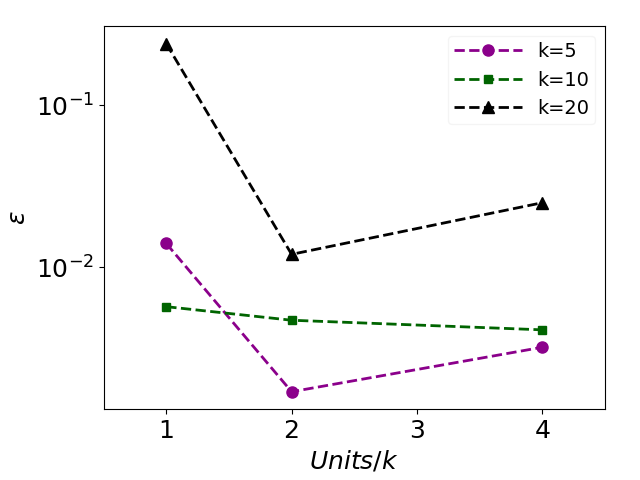}
	\end{minipage}}
	\subfigure[$Layers=2$, PWNN]{
		\begin{minipage}[t]{0.25\textwidth}
			\centering
			\includegraphics[width=1\textwidth]{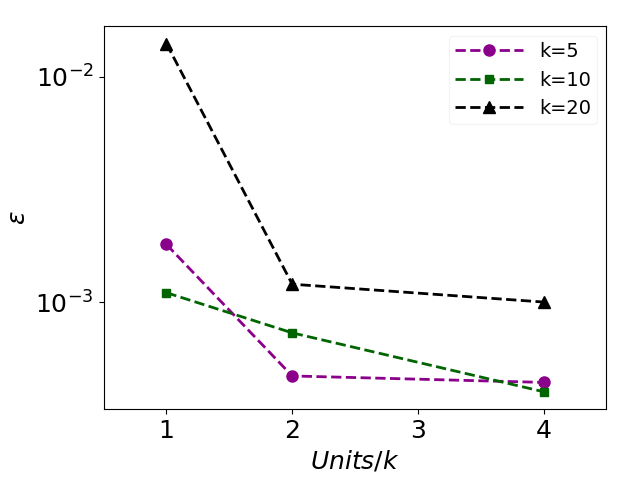}
	\end{minipage}}
	\caption{\textbf{KD problem}. Relative error of TANN, SIREN and PWNN for different $Units$ and wave number $k$, while activation function is fixed in each subfigure. Here, $Layers$ are respectively taken as 1, 2. The x-axis represents $Units/k$. The total numbers of training points are fixed to $N_g=5k\times 4$ and $N_f=5k^2$.}
	\label{error2}
\end{figure*}

\begin{figure*}[htb]
	\centering
	\subfigure[$exact$]{
		\begin{minipage}[t]{0.20\textwidth}
			\centering
			\includegraphics[width=1\textwidth]{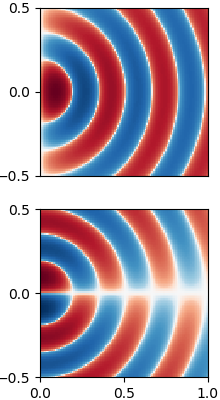}
	\end{minipage}}
	\subfigure[$TANN$]{
		\begin{minipage}[t]{0.20\textwidth}
			\centering
			\includegraphics[width=1\textwidth]{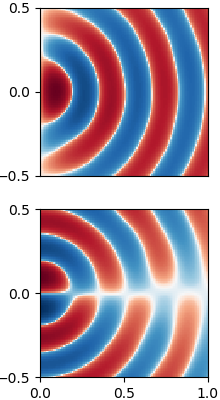}
	\end{minipage}}
	\subfigure[$SIREN$]{
		\begin{minipage}[t]{0.20\textwidth}
			\centering
			\includegraphics[width=1\textwidth]{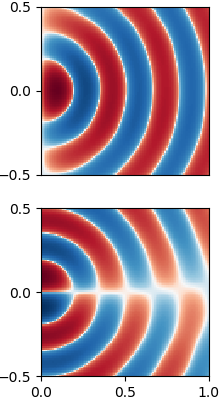}
	\end{minipage}}
	\subfigure[$PWNN$]{
		\begin{minipage}[t]{0.282\textwidth}
			\centering
			\includegraphics[width=1\textwidth]{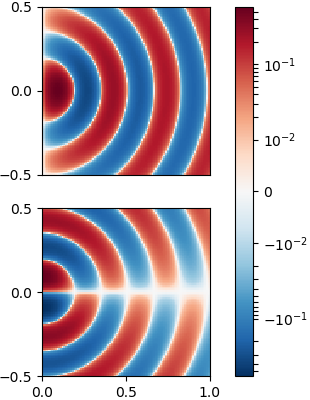}
	\end{minipage}}
	\centering
	\caption{\textbf{KD problem}. The real part(upper part) and imaginary part(lower part) in the whole domain for the exact solution, the solution of TANN, SIREN and PWNN. $Layers=2,Units=40$ for TANN and SIREN, while $Layers=1,Units=40$ for PWNN. Here, the total numbers of training data are fixed to $N_g=100\times 4$ and $N_f=2000$. The wave number in exact solution \ref{es1} is taken as $k=20$.}
	\label{PWNN_fig}
\end{figure*}

It can be seen that the performance of single-layer TANN and SIREN is poor, while the performance of single-layer PWNN is better than that of multi-layer PWNN when $Units$ is not very low. Next, we fix $Units$ and respectivhoely choose $Layers=1,2$ for TANN, SIREN and PWNN then compare them under different numbers of training data.

In Table \ref{TANN_data1}, \ref{TANN_data2}, we fixed the architecture $Layers=1,2,\ Units=40$ for TANN, SIREN and PWNN, and give the relative $L^2$ error under different numbers of training data. Here we randomly generate 50 groups of internal sample points, calculate the average value and standard deviation of relative error. In most cases, we can see that the influence of random selection of internal sample points on the error is acceptable. It can be seen that the error of TANN and SIREN is sensitive to the number of sample points, and the error does not necessarily decrease with the increase of the number of sample points, although the overall trend is downward. This makes the selection of sample points more difficult. Single-layer PWNN is insensitive to the number of sample points, while multi-layer PWNN is sensitive to that. With few internal sample points, the standard deviation of relative error of multi-layer PWNN is large. Obviously, this is because multi-layer PWNN is more difficult to optimize than single-layer PWNN. Until a more suitable optimization method is found, single-layer PWNN is a better choice than multi-layer. 
We know that in the standard finite element method, in order to characterize the wave property of the solution in the region, it is necessary to choose DOFs as $\mathcal{O}(k^2)$. It is also necessary to take $N_f=\mathcal{O}(k^2)$(sample points inside) here in TANN, SIREN and multi-layer PWNN. In Table \ref{TANN_data2}, when $N_f=200$, because of the small number of internal sample points, the expression of internal wave propagation is not enough, which makes the error larger for both TANN and SIREN. But for single-layer PWNN, as it adopts the form of plane wave basis function, we can select much less $N_f$ and achieve a good accuracy. This will reduce a lot of computational cost, because the internal sample points need to calculate the second derivative. In fact, in the later example, we use sinlge-layer PWNN to solve a large wave number problem ($k=100$), and we only choose the training data as $N_g=500\times 4$ and $N_f=500$, this is difficult for other networks to do. As for multi-layer PWNN, although it can also represent the form of plane wave basis function, it also needs more internal sample points due to the increased optimization difficulty.

The relative $L^2$ error of TANN, SIREN and PWNN for different $Units$ and wave number $k$ are shown in Figure \ref{error1}, \ref{error2}. $Layers$ are respectively taken as 1, 2 for TANN, SIREN and PWNN. 
We set the x-axis to $Units/k$, in order to observe can the increased wave number be handled by increasing the network structure accordingly.
In Figure \ref{error1}, $k$ is fixed in each subfigure, PWNN is basically more accurate than TANN and SIREN under the same conditions.
In Figure \ref{error2}, activation function is fixed in each subfigure.
It can be seen that with the increase of $k$, the accuracy of single-layer PWNN increases under the same $Units/k$, while that of TANN and SIREN decreases. This proves that PWNN is more scalable than TANN and SIREN for solve Helmholtz equations. Moreover, single-layer PWNN has greater advantages when dealing with high frequency problems.

\begin{figure*}[htb]
	\centering
	\subfigure[$Layers=1$]{
		\begin{minipage}[t]{0.45\textwidth}
			\centering
			\includegraphics[width=1\textwidth]{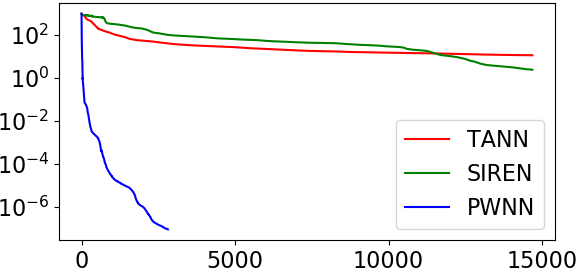}
	\end{minipage}}
	\subfigure[$Layers=2$]{
		\begin{minipage}[t]{0.45\textwidth}
			\centering
			\includegraphics[width=1\textwidth]{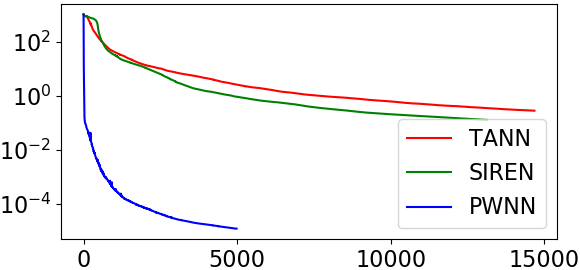}
	\end{minipage}}
	\caption{\textbf{KD problem}. The total loss $\mathcal{M}$ of L-BFGS optimization during the training process. Because the initial loss of different networks is different, we multiply a constant to make the initial loss the same, and observe the descend of loss. The x-axis represents the iteration steps. Here, $Layers$ are respectively taken as 1, 2, and $Units=40$ for TANN, SIREN and PWNN. The wave number in exact solution \ref{es1} is taken as $k=20$. The total numbers of training points are fixed to $N_g=100\times 4$ and $N_f=2000$.}
	\label{loss}
\end{figure*}

Figure \ref{loss} shows the training process of TANN, SIREN and PWNN. $Layers$ are respectively taken as 1, 2, and $Units=40$ for TANN, SIREN and PWNN. It is obvious that the loss of PWNN decreases faster and converges faster than that of TANN and SIREN.
 
Based on the above experimental results, we can draw a conclusion that PWNN works much better than TANN and SIREN on varying architectures or the number of training samples for solving Helmholtz equation.

We now compare one hidden layer PWNN with PWPUM and PWPUM-WT in the domain $\Omega$. 
Here $Units$ indicates the units in PWNN's hidden layer and the number of plane wave basis in PWPUM at the same time. Relative $L^2$ error of PWPUM, PWPUM-WT and PWNN under different $Units$ is shown in Table \ref{unknown_compare10}, \ref{unknown_compare100}, while the numbers of training points in PWNN are maintained as $N_g=50\times 4$ and $N_f=500$ for the wave number $k=10$ in exact solution \ref{es1} and $N_g=500\times 4$ and $N_f=500$ for $k=100$. As expected, more $Units$ corresponds with smaller relative $L^2$ errors for PWPUM, PWPUM-WT and PWNN, PWPUM-WT performs better than PWPUM, and PWNN gives a best output among these three methods under the same number of $Units$. 

\begin{figure}[htbp]
	\centering
	\subfigure[$\varepsilon$ when $k=10$]{
		\begin{minipage}[t]{0.4\textwidth}
			\centering
			\includegraphics[width=1\textwidth]{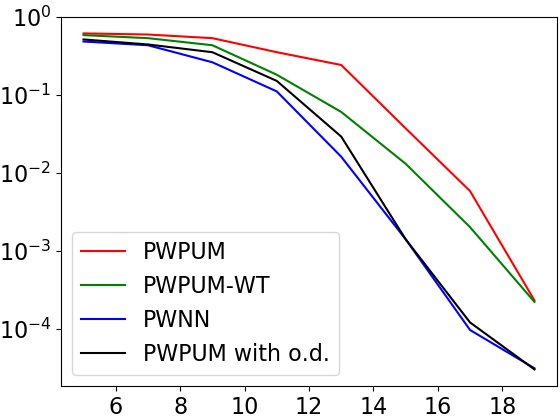}
	\end{minipage}}
	\subfigure[$\varepsilon$ when $k=100$]{
		\begin{minipage}[t]{0.4\textwidth}
			\centering
			\includegraphics[width=1\textwidth]{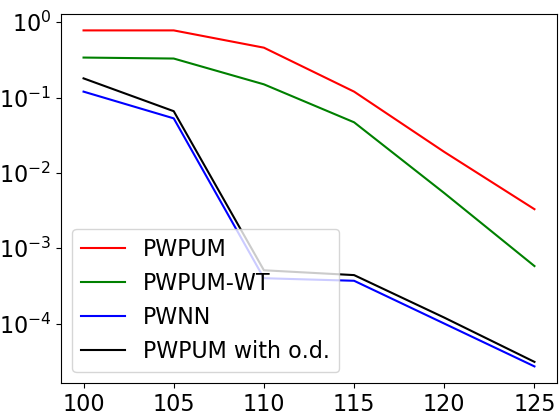}
	\end{minipage}}
	\centering
	\caption{\textbf{KD problem}. Relative $L^2$ error of PWPUM, PWPUM-WT, PWNN and `PWPUM with o.d. (optimized directions)' for different number of $Units$. The x-axis represents the number of $Units$, while y-axis represents $\varepsilon$. The wave number of exact solution \ref{es1} are respectively taken as $k=10$ and $k=100$. Here PWNN has one hidden layer, the total numbers of training data are fixed to $N_g=50\times 4, N_f=500$ for $k=10$ and $N_g=500\times 4, N_f=500$ for $k=100$. }
	\label{two_cases_curve}
\end{figure}

The difference between PWPUM and PWNN is that PWNN will learn a group of optimized plane waves $e^{\mathbf{i} w_i^T \mathbf{x}}$ . After training PWNN with optimized $\{w_i\}$, we can fix the directions $k_i$ in PWPUM \eqref{PWPUM} to $(k{w_{i,1}}/{||w_i||},k{w_{i,2}}/{||w_i||})$. Then use PWPUM to test whether using these directions can get smaller relative $L^2$ errors than the uniformly selected directions, we call this method `PWPUM with o.d. (optimized directions)'. We compute the numerical solutions in two cases of exact solution \ref{es1}: $k=10$ and $k=100$. The relative $L^2$ error of PWPUM, PWPUM-WT, PWNN and `PWPUM with o.d.' for growing number of $Units$ in these two cases is shown in Figure \ref{two_cases_curve}. It can be seen that PWNN has the same convergence order but less relative $L^2$ error than PWPUM and PWPUM-WT. As `PWPUM with o.d.' has higher accuracy than PWPUM and PWPUM-WT, we come to a conclusion that optimized directions are better than the uniformly selected directions or add a rotation angle to it. 
These experimental results can guarantee that PWNN has competitive performance with PWPUM for the KD problem.  

\begin{table}[htbp]
	\centering
	\caption{\textbf{KD problem}. Relative $L^2$ error $\varepsilon$ (Accuracy improvement compared with PWPUM) of PWPUM, PWPUM-WT and PWNN under different $Units$. Here, the total numbers of training data in PWPUM-WT and PWNN are fixed to $N_g=50\times 4$ and $N_f=500$. The wave number in exact solution \ref{es1} is taken as $k=10$.}
	\label{unknown_compare10}
	\small{
	\begin{tabular}{c|ccc}
		\toprule
		$Units$ & PWPUM & PWPUM-WT & PWNN \\
		\midrule
		5 & 6.1e-1 & 5.8e-1($+$0.02) & 4.8e-1($+$0.10) \\
		7 & 5.9e-1 & 5.3e-1($+$0.05) & 4.3e-1($+$0.14) \\
		9 & 5.3e-1 & 4.3e-1($+$0.09) & 2.6e-1($+$0.31) \\
		11 & 3.5e-1 & 1.8e-1($+$0.29) & 1.1e-1($+$0.50) \\
		13 & 2.4e-1 & 6.0e-2($+$0.60) & 1.6e-2($+$1.18) \\
		15 & 3.7e-2 & 1.3e-2($+$0.45) & 1.4e-3($+$1.42) \\
		17 & 5.8e-3 & 2.0e-3($+$0.46) & 9.6e-5($+$1.78) \\
		19 & 2.3e-4 & 2.2e-4($+$0.02) & 3.1e-5($+$0.87) \\
		\bottomrule
	\end{tabular}
	}
\end{table}

\begin{table}[htbp]
	\centering
	\caption{\textbf{KD problem}. Relative $L^2$ error $\varepsilon$ (Accuracy improvement compared with PWPUM) of PWPUM, PWPUM-WT and PWNN under different $Units$. Here, the total numbers of training data in PWPUM-WT and PWNN are fixed to $N_g=500\times 4$ and $N_f=500$. The wave number in exact solution \ref{es1} is taken as $k=100$.}
	\label{unknown_compare100}
	\small{
	\begin{tabular}{c|ccc}
		\toprule
		$Units$ & PWPUM & PWPUM-WT & PWNN \\
		\midrule
		100 & 7.8e-1 & 3.4e-1($+$0.36) & 1.2e-1($+$0.81) \\
		105 & 7.8e-1 & 3.3e-1($+$0.37) & 5.3e-2($+$1.17) \\
		110 & 4.6e-1 & 1.5e-1($+$0.49) & 4.0e-4($+$3.06) \\
		115 & 1.2e-1 & 4.7e-2($+$0.41) & 3.7e-4($+$2.51) \\
		120 & 1.9e-2 & 5.4e-3($+$0.55) & 1.0e-4($+$2.28) \\
		125 & 3.3e-3 & 5.8e-4($+$0.76) & 2.7e-5($+$2.09) \\
		\bottomrule
	\end{tabular}
	}
\end{table}

\subsection{Problem with unknown directions (UD)}
Actually, the exact solution of \eqref{Helm-BV} usually does not cover all directions. We focus a more practical problem, the solution of which only integrate the plane waves with some unknown directions. Consider a square domain $\Omega=[-1,1]\times[-1,1]$, let the exact solution be
\beq\label{es2}
u_*=\sum_{i=1}^d e^{\mathbf{i} k_i^T \mathbf{x}},||k_i||=k,
\eeq
each $k_i$ is a randomly selected two-dimensional vector satisfies $||k_i||=k$, $d$ denotes the number of directions of exact solution, and the coefficient of each direction is 1. 
Our training dataset consists of $N_g=50\times 4$ and $N_f=500$. 

To assess the strength of our approach, we give several different groups $(d,Units)$, $d$ denotes the number of directions of exact solution and $Units$ denotes the units in PWNN's hidden layer and the number of plane wave basis in PWPUM at the same time. We let $k=5,10$ respectively and choose four groups of $(d,Units)$: $(10,10),(10,15),(15,10),(15,15)$ for each $k$.
In each group, we randomly generate 50 exact solutions with $d$ directions, then use PWPUM, PWPUM-WT and PWNN with the corresponding number of $Units$ to solve these corresponding Helmholtz equations. 

\begin{table}[htbp]
	\centering
	\caption{\textbf{UD problem}. Relative $L^2$ errors $\varepsilon$ (Accuracy improvement compared with PWPUM) for four groups of $d,Units$. `Average' represents the average of these relative $L^2$ errors in 50 times. `Max $\varepsilon$ in PWNN' or `Min $\varepsilon$ in PWNN' represents the relative $L^2$ errors between the predicted and the exact solution that makes the PWNN solution's $\varepsilon$ the largest or smallest in this group, respectively. The wave number is taken as $k=5$.}
	\label{compare_unknown,k=5}
	\small{
	\begin{tabular}{c|ccc}
		\toprule
		$d=10,Units=10$ & PWPUM & PWPUM-WT & PWNN \\
		\midrule
		Average & 2.9e-1 & 1.2e-1($+$0.38) & 3.3e-5($+$3.94) \\
		Max $\varepsilon$ in PWNN & 2.6e-1 & 9.3e-2($+$0.45) & 3.9e-4($+$2.82) \\
		Min $\varepsilon$ in PWNN & 4.3e-1 & 1.2e-1($+$0.55) & 1.1e-8($+$7.59) \\
		\midrule
		$d=10,Units=15$ & PWPUM & PWPUM-WT & PWNN \\
		\midrule
		Average & 1.3e-1 & 1.3e-2($+$1.00) & 4.4e-6($+$4.47) \\
		Max $\varepsilon$ in PWNN & 7.6e-2 & 7.7e-3($+$0.99) & 2.7e-5($+$3.45) \\
		Min $\varepsilon$ in PWNN & 4.2e-2 & 1.8e-2($+$0.37) & 2.5e-8($+$6.23) \\
		\midrule
		$d=15,Units=10$ & PWPUM & PWPUM-WT & PWNN \\
		\midrule
		Average & 2.4e-1 & 9.5e-2($+$0.40) & 1.0e-3($+$2.38) \\
		Max $\varepsilon$ in PWNN & 7.0e-2 & 5.2e-2($+$0.13) & 4.3e-3($+$1.21) \\
		Min $\varepsilon$ in PWNN & 3.0e-1 & 1.6e-1($+$0.27) & 3.7e-5($+$3.91) \\
		\midrule
		$d=15,Units=15$ & PWPUM & PWPUM-WT & PWNN \\
		\midrule
		Average & 1.1e-1 & 1.2e-2($+$0.96) & 4.7e-6($+$4.37) \\
		Max $\varepsilon$ in PWNN & 2.8e-2 & 9.5e-3($+$0.47) & 4.2e-5($+$2.82) \\
		Min $\varepsilon$ in PWNN & 3.4e-2 & 1.5e-2($+$0.36) & 9.8e-8($+$5.54) \\
		\bottomrule
	\end{tabular}
	}
\end{table}

\begin{table}[htbp]
	\centering
	\caption{\textbf{UD problem}. 
		Relative $L^2$ errors $\varepsilon$ (Accuracy improvement compared with PWPUM) for four groups of $d,Units$. 
		See the caption of Table \ref{compare_unknown,k=5} for the same notations. 
		The wave number is taken as $k=10$.}
	\label{compare_unknown,k=10}
	\small{
		\begin{tabular}{c|ccc}
			\toprule
			$d=10,Units=10$ & PWPUM & PWPUM-WT & PWNN \\
			\midrule
			Average & 5.8e-1 & 4.8e-1($+$0.08) & 3.0e-1($+$0.29) \\
			Max $\varepsilon$ in PWNN & 6.4e-1 & 4.4e-1($+$0.16) & 4.8e-1($+$0.12) \\
			Min $\varepsilon$ in PWNN & 7.0e-1 & 4.4e-1($+$0.20) & 1.7e-3($+$2.61) \\
			\midrule
			$d=10,Units=15$ & PWPUM & PWPUM-WT & PWNN \\
			\midrule
			Average & 4.4e-1 & 2.8e-1($+$0.20) & 8.8e-4($+$2.70) \\
			Max $\varepsilon$ in PWNN & 3.8e-1 & 3.0e-1($+$0.10) & 2.1e-2($+$0.26) \\
			Min $\varepsilon$ in PWNN & 4.9e-1 & 2.6e-1($+$0.28) & 1.7e-9($+$8.46) \\
			\midrule
			$d=15,Units=10$ & PWPUM & PWPUM-WT & PWNN \\
			\midrule
			Average & 5.9e-1 & 4.9e-1($+$0.08) & 3.0e-1($+$0.29) \\
			Max $\varepsilon$ in PWNN & 6.5e-1 & 4.6e-1($+$0.15) & 6.3e-1($+$0.01) \\
			Min $\varepsilon$ in PWNN & 6.1e-1 & 5.2e-1($+$0.07) & 7.0e-2($+$0.94) \\
			\midrule
			$d=15,Units=15$ & PWPUM & PWPUM-WT & PWNN \\
			\midrule
			Average & 4.0e-1 & 2.6e-1($+$0.20) & 4.3e-3($+$1.97) \\
			Max $\varepsilon$ in PWNN & 4.3e-1 & 2.0e-1($+$0.33) & 2.4e-2($+$1.25) \\
			Min $\varepsilon$ in PWNN & 4.5e-1 & 2.4e-1($+$0.27) & 1.8e-7($+$6.40) \\
			\bottomrule
	     \end{tabular}
    }
\end{table}

The relative $L^2$ errors for each group are summarized in Table \ref{compare_unknown,k=5}, \ref{compare_unknown,k=10}. `Average' represents the average of these relative $L^2$ errors in 50 times. `Max $\varepsilon$ in PWNN' or `Min $\varepsilon$ in PWNN' represents the relative $L^2$ errors between the predicted and the exact solution that makes the PWNN solution's $\varepsilon$ the largest or smallest in this group, respectively.
 
When $k=5$, PWNN works much better than PWPUM and PWPUM-WT in every group. For $k=10$, because $k$ increases and the approximation space remains unchanged, the errors of various methods are increasing. PWNN keeps good accuracy when $Units=15$, while PWPUM and PWPUM-WT can not achieve the relative error less than 0.1. How to use as few $Units$ as possible to learn the direction of the exact solution in the case of large wave number will be our future research work.

\begin{figure}[htbp]
	\centering
	\subfigure[Max,$(d,Units,k)=(10,10,5)$]{
		\begin{minipage}[t]{0.3\textwidth}
			\centering
			\includegraphics[width=1\textwidth]{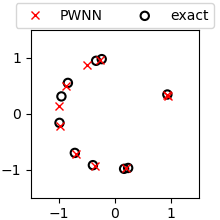}
	\end{minipage}}
	\subfigure[Min,$(d,Units,k)=(10,10,5)$]{
		\begin{minipage}[t]{0.3\textwidth}
			\centering
			\includegraphics[width=1\textwidth]{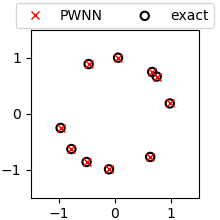}
	\end{minipage}}    
    \subfigure[Max,$(d,Units,k)=(15,10,5)$]{
    	\begin{minipage}[t]{0.3\textwidth}
    		\centering
    		\includegraphics[width=1\textwidth]{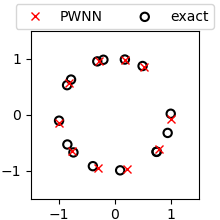}
    \end{minipage}}

    \subfigure[Min,$(d,Units,k)=(15,10,5)$]{
    	\begin{minipage}[t]{0.3\textwidth}
    		\centering
    		\includegraphics[width=1\textwidth]{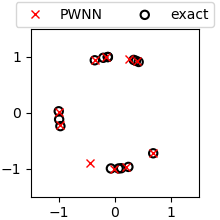}
    \end{minipage}}
    \subfigure[Max,$(d,Units,k)=(10,15,5)$]{
    	\begin{minipage}[t]{0.3\textwidth}
    		\centering
    		\includegraphics[width=1\textwidth]{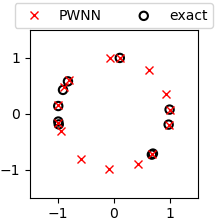}
    \end{minipage}}
    \subfigure[Min,$(d,Units,k)=(10,15,5)$]{
    	\begin{minipage}[t]{0.3\textwidth}
    		\centering
    		\includegraphics[width=1\textwidth]{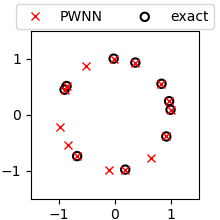}
    \end{minipage}}
  
    \subfigure[Max,$(d,Units,k)=(15,15,5)$]{
    	\begin{minipage}[t]{0.3\textwidth}
    		\centering
    		\includegraphics[width=1\textwidth]{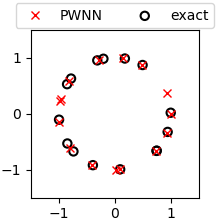}
    \end{minipage}}
    \subfigure[Min,$(d,Units,k)=(15,15,5)$]{
    	\begin{minipage}[t]{0.3\textwidth}
    		\centering
    		\includegraphics[width=1\textwidth]{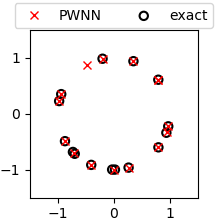}
    \end{minipage}}
    \subfigure[Max,$(d,Units,k)=(10,15,10)$]{
    	\begin{minipage}[t]{0.3\textwidth}
    		\centering
    		\includegraphics[width=1\textwidth]{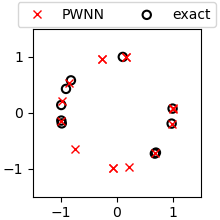}
    \end{minipage}}

    \subfigure[Min,$(d,Units,k)=(10,15,10)$]{
    	\begin{minipage}[t]{0.3\textwidth}
    		\centering
    		\includegraphics[width=1\textwidth]{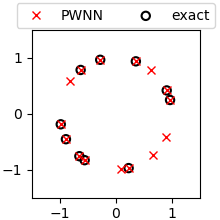}
    \end{minipage}}
    \subfigure[Max,$(d,Units,k)=(15,15,10)$]{
    	\begin{minipage}[t]{0.3\textwidth}
    		\centering
    		\includegraphics[width=1\textwidth]{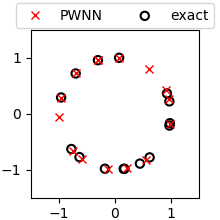}
    \end{minipage}}
    \subfigure[Min,$(d,Units,k)=(15,15,10)$]{
    	\begin{minipage}[t]{0.3\textwidth}
    		\centering
    		\includegraphics[width=1\textwidth]{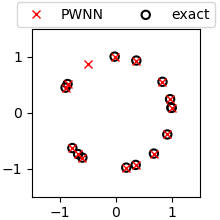}
    \end{minipage}}

	\caption{\textbf{UD problem}. The directions $(\cos\theta,\sin\theta)$ of exact solutions and PWNN solutions in `Max' (`Max $\varepsilon$ in PWNN') and `Min' (`Min $\varepsilon$ in PWNN') cases.}
	\label{direction}
\end{figure}

Figure \ref{direction} prints the directions $(\cos\theta,\sin\theta)$ of exact solutions and PWNN solutions in `Max' (`Max $\varepsilon$ in PWNN') and `Min' (`Min $\varepsilon$ in PWNN') cases of $(d,Units)=(10,10)$,$(15,10)$
,$(10,15)$,$(15,15),k=5$ and $(d,Units)=(10,15),(15,15),k=10$. 
In `Min' case, the optimized directions are almost the same as the directions of exact solution. In `Max' case, the optimized directions also basically the same as most of the directions of exact solution, except some very close directions.

\begin{figure*}[htbp]
	\centering
	\subfigure[$\varepsilon,d=20$]{
		\begin{minipage}[t]{0.45\textwidth}
			\centering
			\includegraphics[width=1\textwidth]{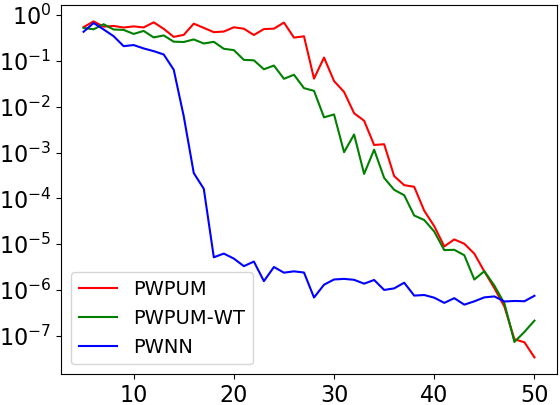}
	\end{minipage}}
	\subfigure[$\varepsilon,d=100$]{
		\begin{minipage}[t]{0.45\textwidth}
			\centering
			\includegraphics[width=1\textwidth]{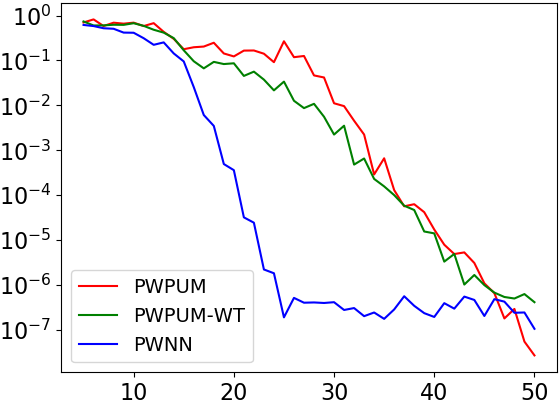}
	\end{minipage}}
	\caption{\textbf{UD problem}. Relative $L^2$ error of PWPUM, PWNN for different number of $Units$. The x-axis represents the number of $Units$, while y-axis represents $\varepsilon$. Let $k=10, d=20$ or $k=10, d=100$ in the exact solutions \eqref{es2}. Here PWNN has one hidden layer, the total numbers of training data are fixed to $N_g=50\times 4$ and $N_f=500$.}
	\label{multiwave_losscurve}
\end{figure*}

Next, we randomly generate exact solutions in two cases: $k=10,d=20$ and $k=10,d=100$. The total numbers of training points in PWNN are fixed to $N_g=50\times 4$ and $N_f=500$. 
$\varepsilon$ of PWPUM, PWPUM-WT and one hidden layer PWNN for growing number of $Units$ is shown in Figure \ref{multiwave_losscurve}. It is obvious that when $Units$ is not very large, e.g. $Units<30$, PWNN gives much better solution than PWPUM. Since it is theoretically guaranteed that the upper bound of relative $L^2$ error of PWPUM will decrease with the increase of $Units$, and PWNN itself has the error of generalization and optimization, PWPUM works better than PWNN when $Units$ is large, e.g. $Units>45$. 
However, to achieve an acceptable relative  $L^2$ error, e.g. 1e-3, 
PWNN only needs to use a much smaller $Units$ than PWPUM.
In UD problem, the performance of PWPUM-WT are not much better than that of PWPUM, which proves that in some cases, adding only one degree of freedom of rotation angle is not enough.
 
Based on the experimental results, we can conclude that, unlike using the plane wave basis with fixed directions in PWPUM, PWNN can find almost the best directions. Thus PWNN is indeed a generalization of PWPUM but more practical. 

\section{Conclusion}
This paper proposed a plane wave activation based neural network (PWNN) for Helmholtz equation. 
We compared PWNN to traditional activation based neural network (TANN), sin based neural network (SIREN) and plane wave partition of unity method (PWPUM). 
We draw lessons from LSM-WT and propose PWPUM-WT, which can be regarded as a simplified PWNN.
The experimental results can guarantee the improvement of introducing plane wave activation for Helmholtz equation and the improvement of PWNN on PWPUM for more practical problems. 
This paper only provide a basic technique for Helmholtz equation. 
We will develop PWNN for some general Helmholtz equations at the future. 

\section*{Acknowledgements}
This research was partially supported by National Natural Science Foundation of China with grant 11831016 and was partially supported by Science Challenge Project, China TZZT2019-B1.1.

\medskip
\small

\bibliographystyle{plain}
\bibliography{yourbibfile}

\begin{thebibliography}{10}

\bibitem{2014A}
M.~Amara, S.~Chaudhry, J.~Diaz, R.~Djellouli, and S.~L. Fiedler.
\newblock A local wave tracking strategy for efficiently solving mid- and
  high-frequency helmholtz problems.
\newblock {\em Computer Methods in Applied Mechanics $\&$ Engineering},
  276:473--508, 2014.

\bibitem{babuska1997the}
Ivo~M {Babuska} and J.~M. {Melenk}.
\newblock The partition of unity method.
\newblock {\em International Journal for Numerical Methods in Engineering},
  40(4):727--758, 1997.

\bibitem{2011Numerical}
Jean-David Benamou, Francis Collino, and Simon Marmorat.
\newblock Numerical microlocal analysis revisited.
\newblock 01 2011.

\bibitem{2004Numerical}
Jean~David Benamou, Francis Collino, and Olof Runborg.
\newblock Numerical microlocal analysis of harmonic wavefields.
\newblock {\em Journal of Computational Physics}, 199(2):717--741, 2004.

\bibitem{berg2018unified}
Jens Berg and Kaj Nystr{\"o}m.
\newblock A unified deep artificial neural network approach to partial
  differential equations in complex geometries.
\newblock {\em Neurocomputing}, 317:28--41, 2018.

\bibitem{buffa2008error}
Annalisa {Buffa} and Peter {Monk}.
\newblock Error estimates for the ultra weak variational formulation of the
  helmholtz equation.
\newblock {\em Mathematical Modelling and Numerical Analysis}, 42(6):925--940,
  2008.

\bibitem{cocquet2017a}
Pierre-Henri {Cocquet}, Martin~J. {Gander}, and Xueshuang {Xiang}.
\newblock A finite difference method with optimized dispersion correction for
  the helmholtz equation.
\newblock {\em International Conference on Domain Decomposition Methods}, pages
  205--213, 2017.

\bibitem{cocquet2019a}
Pierre-Henri {Cocquet}, Martin~J. {Gander}, and Xueshuang {Xiang}.
\newblock Dispersion correction for helmholtz in 1d with piecewise constant
  wavenumber.
\newblock {\em International Conference on Domain Decomposition Methods}, 2019.

\bibitem{cybenko1989approximation}
George Cybenko.
\newblock Approximation by superpositions of a sigmoidal function.
\newblock {\em Mathematics of control, signals and systems}, 2(4):303--314,
  1989.

\bibitem{deraemaeker1999dispersion}
Arnaud {Deraemaeker}, Ivo~M {Babuska}, and Philippe {Bouillard}.
\newblock Dispersion and pollution of the fem solution for the helmholtz
  equation in one, two and three dimensions.
\newblock {\em International Journal for Numerical Methods in Engineering},
  46(4):471--499, 1999.

\bibitem{weinan2017deep}
Weinan E, Jiequn Han, and Arnulf Jentzen.
\newblock Deep learning-based numerical methods for high-dimensional parabolic
  partial differential equations and backward stochastic differential
  equations.
\newblock {\em Communications in Mathematics and Statistics}, 5(4):349--380,
  2017.

\bibitem{2017Learning}
Jun Fang, Jianliang Qian, and Leonardo Zepeda-Nú?ez….
\newblock Learning dominant wave directions for plane wave methods for
  high-frequency helmholtz equations.
\newblock {\em Research in the Mathematical ences}, 4(1), 2017.

\bibitem{gittelson2009plane}
Claude~J. {Gittelson}, Ralf {Hiptmair}, and Ilaria {Perugia}.
\newblock Plane wave discontinuous galerkin methods: Analysis of the h-version.
\newblock {\em Mathematical Modelling and Numerical Analysis}, 43(2):297--331,
  2009.

\bibitem{greydanus2019hamiltonian}
Sam Greydanus, Misko Dzamba, and Jason Yosinski.
\newblock Hamiltonian neural networks.
\newblock {\em arXiv: Neural and Evolutionary Computing}, 2019.

\bibitem{hiptmair2011plane}
R.~{Hiptmair}, A.~{Moiola}, and I.~{Perugia}.
\newblock Plane wave discontinuous galerkin methods for the 2d helmholtz
  equation: Analysis of the $p$-version.
\newblock {\em SIAM Journal on Numerical Analysis}, 49(1):264--284, 2011.

\bibitem{hornik1991approximation}
Kurt Hornik.
\newblock Approximation capabilities of multilayer feedforward networks.
\newblock {\em Neural networks}, 4(2):251--257, 1991.

\bibitem{hornik1989multilayer}
Kurt Hornik, Maxwell Stinchcombe, and Halbert White.
\newblock Multilayer feedforward networks are universal approximators.
\newblock {\em Neural networks}, 2(5):359--366, 1989.

\bibitem{ihlenburg1995finite}
Frank Ihlenburg and Ivo Babuska.
\newblock Finite element solution of the helmholtz equation with high wave
  number part i: The h-version of the fem.
\newblock {\em Computers $\&$ Mathematics With Applications}, 30(9):9--37, 1995.

\bibitem{lecun2012efficient}
Yann~A LeCun, L{\'e}on Bottou, Genevieve~B Orr, and Klaus-Robert M{\"u}ller.
\newblock Efficient backprop.
\newblock In {\em Neural networks: Tricks of the trade}, pages 9--48. Springer,
  2012.

\bibitem{monk1999a}
Peter Monk and Daqing Wang.
\newblock A least-squares method for the helmholtz equation.
\newblock {\em Computer Methods in Applied Mechanics and Engineering},
  175:121--136, 1999.

\bibitem{newman1995frequency-domain}
Gregory~A Newman and David~L Alumbaugh.
\newblock Frequency-domain modelling of airborne electromagnetic responses
  using staggered finite differences.
\newblock {\em Geophysical Prospecting}, 43(8):1021--1042, 1995.

\bibitem{nocedal1980updating}
Jorge Nocedal.
\newblock Updating quasi-newton matrices with limited storage.
\newblock {\em Mathematics of Computation}, 35(151):773--782, 1980.

\bibitem{Partridge1992The}
P.~W. Partridge, C.~A. Brebbia, and L.~C. Wrobel.
\newblock {\em The Dual Reciprocity Boundary Element Method}.
\newblock Elsevier Applied Science, 1992.

\bibitem{Preston1984Finite}
Preston and T.
\newblock Finite elements for electrical engineers.
\newblock {\em Computer Aided Engineering Journal}, 1(5):164, 1984.

\bibitem{raissi2019physics}
M~Raissi, P~Perdikaris, and GE~Karniadakis.
\newblock Physics-informed neural networks: A deep learning framework for
  solving forward and inverse problems involving nonlinear partial differential
  equations.
\newblock {\em Journal of Computational Physics}, 378:686--707, 2019.

\bibitem{rudd2013constrained}
Keith Rudd, Gianluca Di~Muro, and Silvia Ferrari.
\newblock A constrained backpropagation approach for the adaptive solution of
  partial differential equations.
\newblock {\em IEEE transactions on neural networks and learning systems},
  25(3):571--584, 2013.

\bibitem{shen2007an}
Liang Shen and Yijun Liu.
\newblock An adaptive fast multipole boundary element method for
  three-dimensional acoustic wave problems based on the burton-miller
  formulation.
\newblock {\em Computational Mechanics}, 40(3):461--472, 2007.

\bibitem{sirignano2018dgm}
Justin Sirignano and Konstantinos Spiliopoulos.
\newblock Dgm: A deep learning algorithm for solving partial differential
  equations.
\newblock {\em Journal of Computational Physics}, 375:1339--1364, 2018.

\bibitem{2020Implicit}
Vincent Sitzmann, Julien N.~P Martel, Alexander~W Bergman, David~B Lindell, and
  Gordon Wetzstein.
\newblock Implicit neural representations with periodic activation functions.
\newblock 2020.

\bibitem{trabelsi2018deep}
Chiheb {Trabelsi}, Olexa {Bilaniuk}, Ying {Zhang}, Dmitriy {Serdyuk}, Sandeep
  {Subramanian}, Joao~Felipe {Santos}, Soroush {Mehri}, Negar {Rostamzadeh},
  Yoshua {Bengio}, and Christopher~J {Pal}.
\newblock Deep complex networks.
\newblock In {\em ICLR 2018 : International Conference on Learning
  Representations 2018}, 2018.

\bibitem{yang2020b-pinns:}
Liu Yang, Xuhui Meng, and George~Em Karniadakis.
\newblock B-pinns: Bayesian physics-informed neural networks for forward and
  inverse pde problems with noisy data.
\newblock {\em arXiv: Machine Learning}, 2020.

\bibitem{zang2019weak}
Yaohua Zang, Gang Bao, Xiaojing Ye, and Haomin Zhou.
\newblock Weak adversarial networks for high-dimensional partial differential
  equations.
\newblock {\em arXiv: Numerical Analysis}, 2019.

\end{thebibliography}

\end{document}